\documentclass[10pt,a4paper]{amsart}
\usepackage[applemac]{inputenc} 
\usepackage{latexsym}
\usepackage{color,graphicx,shortvrb}
\usepackage{amsmath, amssymb}
\usepackage{amsfonts}
\usepackage[colorlinks, bookmarks=true]{hyperref}
\newtheorem{theorem}{Theorem}[section]
\newtheorem{lemma}[theorem]{Lemma}

\theoremstyle{definition}

\author{J. M. Almira and  Kh. F. Abu-Helaiel}
\title{ A qualitative description of graphs of discontinuous polynomial functions}

\begin{document}




\begin{abstract} We prove that, if  $f:\mathbb{R}^n\to\mathbb{R}$ satisfies  Fr\'{e}chet's functional equation
\[
\Delta_h^{m+1}f(x)=0 \ \ \text{ for all  }x=(x_1,\cdots,x_n),h=(h_1,\cdots,h_n) \in\mathbb{R}^n,
\]
and $f(x_1,\cdots,x_n)$ is not an ordinary algebraic polynomial in the variables $x_1,\cdots,x_n$, then   $f$ is unbounded on all non-empty open set $U\subseteq \mathbb{R}^n$. Furthermore,  the set $\overline{G(f)}^{\mathbb{R}^{n+1}}$ contains an unbounded open set.  
\end{abstract}

\maketitle

\markboth{J. M. Almira and Kh. F. Abu-Helaiel}{A qualitative description of graphs of discontinuous polynomial functions}
\section{Motivation}

One of the best known functional equations that exists in the literature is Fr\'{e}chet's functional equation, which is given by 
\begin{equation}\label{frechet}
\Delta^{m+1}_hf(x)=0 \ \ (x,h\in X),                                                                                      
\end{equation}
where $f:X\to Y$ denotes a function, $X,Y$ are two $\mathbb{Q}$-vector spaces, and $\Delta^{k}_hf(x)$ is defined inductively by $\Delta_h^1f(x)=f(x+h)-f(x)$ and $\Delta_{h}^{k+1}f(x)=\Delta_{h}^1\left(\Delta_{h}^kf\right)(x)$, $k=1,2,\cdots$. A simple induction argument shows that \eqref{frechet} can be explicitly written as
\begin{equation}\label{frepasofijo}
\Delta_{h}^{m+1}f(x):=\sum_{k=0}^{m+1}\binom{m+1}{k}(-1)^{m+1-k}f(x+kh)=0 \ \ (x,h\in X).
\end{equation}

This  equation was  introduced in the literature by M. Fr\'{e}chet in 1909, for $X=Y=\mathbb{R}$,  as a particular case of the functional equation  
\begin{equation}\label{fre}
\Delta_{h_1h_2\cdots h_{m+1}}f(x)=0 \ \ (x,h_1,h_2,\dots,h_{m+1}\in \mathbb{R}),
\end{equation}
where $f:\mathbb{R}\to\mathbb{R}$ and $\Delta_{h_1h_2\cdots h_s}f(x)=\Delta_{h_1}\left(\Delta_{h_2\cdots h_s}f\right)(x)$, $s=2,3,\cdots$. Indeed, thanks to a classical result by Djokovi\'{c} \cite{Dj}, the equation with variable steps 
$\Delta_{h_1h_2\cdots h_{m+1}}f(x)=0$ is equivalent to the equation with fixed step $\Delta^{m+1}_hf(x)=0$. After Fr\'{e}chet's 
seminal paper \cite{frechet}, the solutions of \eqref{frechet} are named ``polynomial functions'' by the Functional Equations community, since it is known that, under very mild regularity conditions on $f$, if $f:\mathbb{R}\to\mathbb{R}$ satisfies \eqref{frechet}, then $f(x)=a_0+a_1x+\cdots a_{m}x^{m}$ for all $x\in\mathbb{R}$ and certain constants $a_i\in\mathbb{R}$. Indeed, it is known that if $f$ is a solution of $(\ref{frechet})$ with $X=Y=\mathbb{R}$, then $f$ is an ordinary polynomial of degree $\leq m$, $f(x)=a_0+a_1x+\cdots+a_mx^m$, if and only if $f$ is bounded on some set $A\subset\mathbb{R}$  with positive Lebesgue measure $|A|>0$. In particular, all measurable polynomial functions $f:\mathbb{R}\to\mathbb{R}$ are ordinary polynomials. This result was firstly proved for the Cauchy functional equation by Kormes in 1926 \cite{kormes}. Later on, in 1959, the result was proved for polynomials by Ciesielski \cite{ciesielski} (see also \cite{kuczma1}, \cite{kurepa}, \cite{mckiernan}, \cite{laszlo1}). A weaker result is the so called Darboux type theorem, which claims that the polynomial function  $f:\mathbb{R}\to\mathbb{R}$ is an ordinary polynomial if and only if $f_{|(a,b)}$ is  bounded for some nonempty open interval $(a,b)$ (see \cite{darboux}, \cite{sanjuan} for the original result, which was stated for solutions of the Cauchy functional equation and \cite{almira_antonio}, \cite{AK_MJM},  \cite{laszlo1} for a direct proof of this result with polynomial functions).

In \cite{almira_antonio}, \cite{AK_MJM} Fr\'{e}chet's equation was studied from a new fresh perspective. The main idea was to use the basic properties of Lagrange interpolation polynomials in one real variable. This allowed the authors to give a description of the closure of the graph $G(f)=\{(x,f(x)):x\in\mathbb{R}\}$ of any discontinuous polynomial function$f:\mathbb{R}\to\mathbb{R}$. Concretely, they proved that
\begin{equation}\label{grafo}
\overline{G(f)}^{\mathbb{R}^2}=C(l,u)=\{(x,y)\in\mathbb{R}^2: l(x)\leq y\leq u(x)\}
\end{equation}
for a certain  pair of functions  $l,u:\mathbb{R}\to\mathbb{R}\cup\{+\infty,-\infty\}$ such that
\begin{itemize}
\item[(i)] $u$ is lower semicontinuous and $l$ is upper semicontinuous.
\item[(ii)] For all $x\in \mathbb{R}$ we have that $u(x)-l(x)=+\infty$.
\item[(iii)] There exist two non-zero ordinary polynomials $p,q$ such that $p\neq q$ and for all $x\in\mathbb{R}$, we have that  $\{x\}\times [p(x),q(x)]\subseteq C(l,u)$.                                           
\end{itemize}
Clearly, this result implies the Darboux type theorem for the Fr\'{e}chet functional equation. Furthermore, it states that, for every discontinuous polynomial function$f:\mathbb{R}\to\mathbb{R}$, the set $\overline{G(f)}^{\mathbb{R}^2}$ contains an unbounded open set. This is a nice property which stands up, in a very visual form, the fact that discontinuous polynomial function functions have wild oscillations. In this paper we present a new proof of this result, based on the standard  tensor product technique for the Lagrange interpolation problem in several variables, and we use the new focus to prove that, for every $n>1$, if $f:\mathbb{R}^n\to\mathbb{R}$ is a discontinuous polynomial function, then $f$ is locally unbounded and  the closure of its its graph, $\overline{G(f)}^{\mathbb{R}^{n+1}}=\overline{\{(x,f(x):x\in\mathbb{R}^n\}}^{\mathbb{R}^{n+1}}$, contains an unbounded open set. 

Along this paper,  $\Pi_{m,max}^n$ denotes the set of algebraic polynomials in the $n$ variables $x_1,x_2,\cdots,x_n$ with degree $\leq m$ in each one of these variables,
\[
\Pi_{m,max}^n=\{\sum_{0\leq i_1,i_2,\cdots,i_n\leq m}a_{i_1,i_2,\cdots,i_n}x_1^{i_1}x_2^{i_2}\cdots x_{n}^{i_n}:a_{i_1,i_2,\cdots,i_n}\in\mathbb{R}\text{ for all } (i_1,i_2,\cdots,i_n)\}.
\] 
When $n=1$ we simply write $\Pi_m$. 

\section{Main results}

Let  $f:\mathbb{R}^n\to\mathbb{R}$ be an arbitrary function. Take $a=(a_1,\cdots,a_n)\in\mathbb{R}^n$, $h_1,\cdots, h_{n+1}\in\mathbb{R}\setminus\{0\}$, and $\gamma=\{v_k\}_{k=1}^{n+1}\subset \mathbb{R}^n\setminus \{(0,0,\cdots,0)\}$. Then, by tensor product interpolation, it is known that there exists a unique  algebraic polynomial $P(t_1,\cdots,t_{n+1})\in \Pi_{m,\max}^{n+1}$ such that 
\[
P(i_1h_1,i_2h_2,\cdots,i_{n+1}h_{n+1})=f_{i_1,\cdots,i_{n+1}}:=f(a+\sum_{k=1}^{n+1}i_kh_kv_k), \text{ for all } 0\leq i_k\leq m,\ 1\leq k\leq n+1.
\] 
In all what follows, we denote this polynomial  by  $P_{a,h,\gamma}$, where $h:=(h_1,\cdots,h_{n+1})$.  

 
\begin{lemma} \label{MP_lem1} If $f:\mathbb{R}^n\to\mathbb{R}$ satisfies Fr\'{e}chet's functional equation of order $m+1$, 
$\Delta_h^{m+1}f(x)=0 $  for all $x,h \in\mathbb{R}^n$,  then  
\[
P_{a,h,\gamma}(i_1h_1,i_2h_2,\cdots,i_{n+1}h_{n+1})=f(a+\sum_{k=1}^{n+1}i_kh_kv_k), \text{ for all  } (i_1,\cdots,i_{n+1})\in\mathbb{Z}^{n+1}.
\]
\end{lemma}

\noindent \textbf{Proof. } Let us fix the values of  $k\in\{1,\cdots,n+1\}$ and $i_1,i_2,\cdots,i_{k-1},i_{k+1},\cdots,i_{n+1}\in\{0,1,\cdots,m\}$, and let us consider the polynomial of one variable  
$$q_k(x)=P_{a,h,\gamma}(i_1h_1,\cdots, i_{k-1}h_{k-1},x,i_{k+1}h_{k+1}, \cdots,i_{n+1}h_{n+1}) .$$ 
Obviously $q_k\in\Pi_m^1$, so that
\begin{eqnarray*}
0 &=& \Delta_{h_k}^{m+1}q_k(0)=\sum_{r=0}^{m+1}\binom{m+1}{r}(-1)^{m+1-r}q_k(rh_k) \\
&=& \sum_{r=0}^{m}\binom{m+1}{r}(-1)^{m+1-r}P_{a,h,\gamma}(i_1h_1,\cdots, i_{k-1}h_{k-1},rh_k,i_{k+1}h_{k+1}, \cdots,i_{n+1}h_{n+1}) \\
&\ & \ \ \ + q_k((m+1)h_k)\\
&=& \sum_{r=0}^{m}\binom{m+1}{r}(-1)^{m+1-r}f(a+\sum_{(0\leq j\leq n+1;\ j\neq k)} i_jh_jv_j+ rh_kv_k) + q_k((m+1)h_k)\\
&=& \Delta_{h_kv_k}^{m+1}f(a+\sum_{(0\leq j\leq n+1;\ j\neq k)} i_jh_jv_j) -f(a+\sum_{(0\leq j\leq n+1;\ j\neq k)} i_jh_jv_j+ (m+1)h_kv_k)\\
&\ & \ \ \ + q_k((m+1)h_k)\\
&=&  q_k((m+1)h_k)-f(a+\sum_{(0\leq j\leq n+1;\ j\neq k)} i_jh_jv_j+ (m+1)h_kv_k).
\end{eqnarray*}
It follows that  
\begin{equation} \label{nuevodato}
 \begin{array}{cccccccc}
q_k((m+1)h_k) & = & \  &   P_{a,h,\gamma}(i_1h_1,\cdots, i_{k-1}h_{k-1},(m+1)h_k,i_{k+1}h_{k+1}, \cdots,i_{n+1}h_{n+1}) \\
\  & = & \  &  f(a+\sum_{(0\leq j\leq n+1;\ j\neq k)} i_jh_jv_j+ (m+1)h_kv_k).
\end{array}
\end{equation}
Let us now consider the unique  polynomial $P\in\Pi_{m,max}^{n+1}$ which satisfies the Lagrange interpolation conditions 
\[
\begin{array}{llllll}
P(i_1h_1,i_2h_2,\cdots,i_{n+1}h_{n+1}) &=& f(a+\sum_{k=1}^{n+1}i_kh_kv_k), & \text{ for all  } &   0\leq i_j\leq m, \ 1\leq j\leq n+1, j\neq k,\\
\ & \ & \ & \text{ and all } & 1\leq i_k\leq m+1,
\end{array}
\] 
We have already demonstrated, with formula \eqref{nuevodato}, that this polynomial  coincides with  $P_{a,h,\gamma}$. Furthermore, the very same arguments used to prove  \eqref{nuevodato}, applied to the polynomial $P=P_{a,h,\gamma}$, lead us to the conclusion that 
\[
 P_{a,h,\gamma}(i_1h_1,\cdots, i_{k-1}h_{k-1},(m+2)h_k,i_{k+1}h_{k+1}, \cdots,i_{n+1}h_{n+1}) =  f(a+\sum_{(0\leq j\leq n+1;\ j\neq k)} i_jh_jv_j+ (m+2)h_kv_k)
\]
In an analogous way, clearing this time the first term of the sum, and taking as starting point the equality 
\[
 \Delta_{h_kv_k}^{m+1}f(a+\sum_{(0\leq j\leq n+1;\ j\neq k)} i_jh_jv_j-h_kv_k)=0, 
\] 
we conclude that 
\[
 P_{a,h,\gamma}(i_1h_1,\cdots, i_{k-1}h_{k-1},-h_k,i_{k+1}h_{k+1}, \cdots,i_{n+1}h_{n+1}) = f(a+\sum_{(0\leq j\leq n+1;\ j\neq k)} i_jh_jv_j- h_kv_k).
\]
Repeating these arguments forward and backward infinitely many times, and for each  $k\in\{1,\cdots,n+1\}$, we get  
\[
P_{a,h,\gamma}(i_1h_1,i_2h_2,\cdots,i_{n+1}h_{n+1})=f(a+\sum_{k=1}^{n+1}i_kh_kv_k), \text{ for all } (i_1,\cdots,i_{n+1})\in\mathbb{Z}^{n+1},
\] 
which is what we wanted to prove. {\hfill $\Box$}


\begin{lemma}\label{MP_lem2} If $f:\mathbb{R}^n\to\mathbb{R}$ satisfies Fr\'{e}chet's functional equation of order $m+1$, then 
\[
P_{a,h,\gamma}(r_1h_1,r_2h_2,\cdots,r_{n+1}h_{n+1})=f(a+\sum_{k=1}^{n+1}r_kh_kv_k), \text{ for all } (r_1,\cdots,r_{n+1})\in\mathbb{Q}^{n+1}.
\]
Consequently, $\overline{G(f)}^{\mathbb{R}^{n+1}}$ contains the set $\varphi_{\gamma}(\mathbb{R}^{n+1})$, where 
\[
\varphi_{\gamma}(t_1,\cdots,t_{n+1})=(a+\sum_{k=1}^{n+1}t_kv_k,P_{a,h,\gamma}(t_1,\cdots,t_{n+1})). 
\]
\end{lemma}
\noindent \textbf{Proof. } It is enough to take into account that, if $p_1,p_2,\cdots,p_{n+1}\in \mathbb{Z}\setminus\{0\}$, and we use Lemma \ref{MP_lem1}  with the polynomial $P^*(t_1,\cdots,t_{n+1})$ which satisfies the interpolation conditions 
\[
P^*(i_1h_1^*,i_2h_2^*,\cdots,i_{n+1}h_{n+1}^*)=f(a+\sum_{k=1}^{n+1}i_kh_k^*v_k), \text{ for all } 0\leq i_k\leq m,\ 1\leq k\leq n+1,
\] 
where $h_i^*=h_i/p_i$, $i=1,\cdots,n+1$, then 
\begin{eqnarray*}
P^*(i_1h_1,i_2h_2,\cdots,i_{n+1}h_{n+1}) &=& P^*(p_1i_1h_1^*,p_2i_2h_2^*,\cdots,p_{n+1}i_{n+1}h_{n+1}^*) \\
&=& f(a+\sum_{k=1}^{n+1}p_ki_kh_k^*v_k)\\ 
&=& f(a+\sum_{k=1}^{n+1}i_kh_kv_k)\\
&=& P_{a,h,\gamma}(i_1h_1,i_2h_2,\cdots,i_{n+1}h_{n+1}), 
\end{eqnarray*} 
for all  $0\leq i_k\leq m$ and $1\leq k\leq n+1$. Thus, $P^*=P_{a,h,\gamma}$, which implies the first claim in the lemma, since  $p_1,p_2,\cdots,p_{n+1}\in \mathbb{Z}\setminus\{0\}$ were arbitrary. Second claim is a direct consequence of the density of $\mathbb{Q}$ in the real line $\mathbb{R}$.  {\hfill $\Box$}

\begin{lemma}\label{lema2}
Every polynomial $P(x,y)\in\Pi_{m,max}^2$ can be decomposed as
$$P(x,y)=\sum_{i=0}^{2m}A_i(x+y)x^i,$$
where $A_i(t)\in \Pi_m$ is a  polynomial in one variable of degree $\leq m$, for all $i=0,1,\cdots, 2m$.
\end{lemma}

\noindent \textbf{Proof. } Let us consider the change of variables $\varphi(x,y)=(x,x+y)$. If we denote $f_1=x$, $f_2=x+y$, then $y=f_2-f_1$ and, consequently, a simple computation shows that every polynomial 
$P(x,y)=\sum_{i=0}^m\sum_{j=0}^ma_{i,j}x^iy^j\in \Pi_m^2$ can be decomposed as follows:
\begin{eqnarray*}
P(x,y) &= & \sum_{i=0}^m\sum_{j=0}^ma_{i,j}x^iy^j = \sum_{i=0}^m\sum_{j=0}^ma_{i,j}f_1^i(f_2-f_1)^j \\
&=& \sum_{i=0}^m\sum_{j=0}^ma_{i,j}f_1^i (\sum_{s=0}^j\binom{j}{s}(-1)^{j-s}f_2^sf_1^{j-s}) \\
&=& \sum_{i=0}^m\sum_{j=0}^m(\sum_{s=0}^ja_{i,j}\binom{j}{s}(-1)^{j-s}f_2^sf_1^{i+j-s}) \\
&=& \sum_{i=0}^{2m}A_i(f_2)f_1^i \\
&=& \sum_{i=0}^{2m}A_i(x+y)x^i,
\end{eqnarray*}
where $A_i(t)$ is a  polynomial in one variable of degree $\leq m$, for all $i=0,1,\cdots, 2m$.  {\hfill $\Box$}

\begin{theorem}[Description of $G(f)$ for the univariate setting]  \label{MP_grafo_1v} If $f:\mathbb{R}\to\mathbb{R}$ satisties Fr\'{e}chet's functional equation
\[
\Delta_h^{m+1}f(x)=0 \ \ \text{ for all  }(x,h) \in\mathbb{R}^2,
\]
and $f(x)$ is not an ordinary algebraic polynomial, then $f$ is locally unbounded. Indeed, for each  
 $x\in\mathbb{R}$ there exists an unbounded interval  $I_x\subseteq \mathbb{R}$ such that  $\{x\}\times I_x\subseteq  \overline{G(f)}^{\mathbb{R}^2}$. Furthermore, $ \overline{G(f)}^{\mathbb{R}^2}$ contains an unbounded open set. 
\end{theorem}

\noindent \textbf{Proof. }  Let $f:\mathbb{R}\to\mathbb{R}$ be a solution of Fr\'{e}chet's equation  $\Delta_h^{m+1}f=0$. Then Lemma \ref{MP_lem2} guarantees that there exists a unique polynomial  $p_{x_0,h_1,h_2}(x,y)\in  \Pi_{m,max}^2 $ satisfying  
\begin{equation}\label{Q2}
f(x_0+ih_1+jh_2)=p(ih_1,jh_2) \text{ for all }(i,j)\in\mathbb{Q}^2.
\end{equation}
and
 \begin{equation}\label{MP_grafo_1}
 \Gamma_{x_0,h_1,h_2}:=\{(x_0+u+v,p_{x_0,h_1,h_2}(u,v)):u,v\in\mathbb{R}\}\subseteq \overline{G(f)}^{\mathbb{R}^2}, 
 \end{equation}
Thus, we are now interested on studying the sets  $\Gamma_{x_0,h_1,h_2}$. 
If $p_{x_0,h_1,h_2}(x,y)=A(x+y)$ for a certain univariate polynomial $A$, then $\Gamma_{x_0,h_1,h_2}$ has empty interior and, in fact, it coincides with the graph of an ordinary algebraic polynomial. Hence, in this case  the property (\ref{MP_grafo_1}) does not add any extra interesting information. 

We claim that, if  $f$ is not an ordinary algebraic polynomial, there exist $x_0,h_1,h_2\in\mathbb{R}$ such that  $p_{x_0,h_1,h_2}(x,y)$ is not a polynomial in the variable $x+y$. Concretely, we will prove that, for adequate values $x_0,h_1$ and $h_2$, this polynomial admits a decomposition of the form
\begin{equation}\label{MP_poli_estricto}
p_{x_0,h_1,h_2}(x,y)=\sum_{i=0}^NA_i(x+y)x^i,\ \text{ with } A_N(t)\neq 0 \text{ and } N\geq 1,
\end{equation}
where $A_i(t)$ is an univariate polynomial of degree $\leq m$, for  $i=0,1,\cdots,N$.  

Obviously, if  $\eqref{MP_poli_estricto}$ holds true, then for every $\alpha\in\mathbb{R}\setminus Z(A_N)$ (where $Z(A_N)=\{s\in\mathbb{R}:A_N(s)=0\}$ is a finite set with at most $m$ points), we have that  $m_{\alpha}(x)=p_{x_0,h_1,h_2}(x,\alpha-x)=\sum_{i=0}^NA_i(\alpha)x^i$ is a non-constant polynomial, so that  $m_{\alpha}(\mathbb{R})$ is an unbounded interval.  Furthermore,
\[
\{x_0+\alpha\}\times m_{\alpha}(\mathbb{R})\subseteq \Gamma_{x_0,h_1,h_2}\subseteq \overline{G(f)}^{\mathbb{R}^2}.
\]
Thus, if  $\eqref{MP_poli_estricto}$ is satisfied, then  $f$ is locally unbounded and, for each  $x\in\mathbb{R}$ there exists an unbounded interval $I_x\subseteq \mathbb{R}$ such that 
$\{x\}\times I_x\subseteq  \overline{G(f)}^{\mathbb{R}^2}$. 

Let us demonstrate that, if $P=p_{x_0,h_1,h_2}$ satisfies  $\eqref{MP_poli_estricto}$, then  $\Gamma_{x_0,h_1,h_2}$ contains an unbounded  open set. To prove this, we consider the function  $\varphi:\mathbb{R}^2\to\mathbb{R}^2$, 
$$\varphi(x,y)=(x+y+x_0,P(x,y)).$$ 
A simple computation reveals that   $$\det \varphi'(x,y)=P_y-P_x=-\sum_{k=1}^NkA_k(x+y)x^{k-1}$$ is a nonzero polynomial, so that 
 $\Omega=\mathbb{R}^2\setminus \{(x,y): \det \varphi'(x,y)=0\}$ is a non-empty open subset of the plane. Indeed, $\Omega$ is a dense open subset of $\mathbb{R}^2$. Thus, we can apply the Open Mapping Theorem for differentiable functions defined over finite dimensional Euclidean spaces, to the function  $\varphi$, concluding that  $W=\varphi(\Omega)$ is an open subset of  $\mathbb{R}^2$ which is contained into $\Gamma_{x_0,h_1,h_2}$. Furthermore, the inclusions  $\{x_0+\alpha\}\times m_{\alpha}(\mathbb{R})\subseteq \Gamma_{x_0,h_1,h_2}$ prove that $W$ is unbounded. 

Let us now show that, in fact, the relation $\eqref{MP_poli_estricto}$ holds true for certain values $x_0,h_1,h_2$. It follows from Lemma \ref{lema2} that  $p_{x_0,h_1,h_2}(x,y)$ admits a decomposition of the form 
\begin{equation}\label{MP_poli_estricto_nuevo}
p_{x_0,h_1,h_2}(x,y)=\sum_{i=0}^NA_i(x+y)x^i,\ \text{ with } A_N(t)\neq 0 \text{ and } N\geq 0.
\end{equation}
Thus, our claim is that, for certain choice of $x_0,h_1,h_2$, the decomposition \eqref{MP_poli_estricto_nuevo} satisfies $N\geq 1$. Assume, on the contrary, that $N=0$ for all  $x_0,h_1,h_2$. Then, for any fixed pair of values   $h_1,h_2$, every polynomial   $p_{x_0,h_1,h_2}(x,y)$ satisfies a relation of the form $p_{x_0,h_1,h_2}(x,y)=A_{x_0}(x+y)$ for certain polynomial $A_{x_0}\in\Pi_m$. Hence, the assumption that $f$ is not an ordinary algebraic polynomial, implies that there exist two distinct points $x_0,x_1\in\mathbb{R}$ such that  $A_{x_1}(0)\neq A_{x_0}(x_1-x_0)$, since otherwise, if we fix the value $x_0$ and take $x\in\mathbb{R}$ arbitrary, we would have that  
\[
f(x)=p_{x,h_1,h_2}(0)=A_x(0)=A_{x_0}(x-x_0),
\] 
and $f$ would be an ordinary polynomial.

Let us now consider the polynomial  
 $p_{x_0,x_1-x_0,h_2}(x,y)$. By hypothesis, this polynomial satisfies the identity $p_{x_0,x_1-x_0,h_2}(x,y)=A(x+y)$ for certain $A\in \Pi_m$. Now, a simple computation shows that 
\[
A(x_1-x_0)=p_{x_0,x_1-x_0,h_2}(x_1-x_0,0)=f(x_0+(x_1-x_0))=f(x_1)=A_{x_1}(0).
\]
On the other hand, for each $j\in\mathbb{Z}$, we have that 
\[
A(jh_2)=p_{x_0,x_1-x_0,h_2}(0,jh_2)=f(x_0+jh_2)=A_{x_0}(jh_2),
\]
so that $A$ and $A_{x_0}$ coincide in infinitely many points. Thus they are the same polynomial, and  $A_{x_1}(0)= A_{x_0}(x_1-x_0)$, which contradicts the assumption that $f$ is not a polynomial. 
{\hfill $\Box$}

Now we state and prove the main result of this paper:

\begin{theorem}[Description of  $G(f)$ for the multivariate setting]  \label{MP_grafo_vv} If $f:\mathbb{R}^n\to\mathbb{R}$ satisfies Fr\'{e}chet's functional equation 
\[
\Delta_h^{m+1}f(x)=0 \ \ \text{ for all  }x,h \in\mathbb{R}^n,
\]
and $f(x_1,\cdots,x_n)$ is not an ordinary algebraic polynomial, then  $f$ is locally  unbounded. Furthermore,  $\overline{G(f)}^{\mathbb{R}^{n+1}}$ contains an unbounded open set. 
\end{theorem}

Previous to give the formal proof of Theorem  \ref{MP_grafo_vv}, some remarks are necessary. Thus, just to start, we observe that if   $f(x_1,\cdots,x_n)$ is not an ordinary algebraic polynomial, then there exist some values (which we fix from now on) $s\in \{1,\cdots,n\}$ and  $(a_1,a_2,\cdots,a_{s-1},a_{s+1},\cdots,a_n)\in\mathbb{R}^{n-1}$ such that 
$$g(x)=f(a_1,a_2,\cdots,a_{s-1},x,a_{s+1},\cdots,a_n)$$ is not an ordinary algebraic polynomial. This result has been proved in several ways and can be found, for example, in   \cite{AK_CJM},  \cite{kuczma}, and  \cite{prager}. Furthermore, if we take into account the proof of Theorem \ref{MP_grafo_1v}, we know that, if we denote by  $p_{x_0,\alpha,\beta}(x,y)$ 
the unique element of $\Pi_{m,\max}^2$ such that 
\[
p_{x_0,\alpha,\beta}(i\alpha,j\beta)=g(x_0+i\alpha+j\beta),\ \ \text{ for all } i,j=0,1,\cdots,m,
\]
then there exist $a_s,h_s,h_{n+1}\in\mathbb{R}$, $1\leq N\leq 2m$, and polynomials $A_k\in \Pi_m$, $k=0,1,\cdots, N$ such that 
\[
p_{a_s,h_s,h_{n+1}}(x,y)=\sum_{k=0}^NA_k(x+y)x^k,\ \ \text{ and } A_N\neq 0.
\]
We also fix, from now on, the values $a_s,h_s$ and ,$h_{n+1}$. Furthermore, we also fix the values $h_1$, $\cdots$, $h_{s-1}$, $h_{s+1}$, $ \cdots$, $h_n$ with the only imposition that they are all real numbers different from zero.

\begin{lemma}\label{MP_lem3} Let us use, with the values $a=(a_1,\cdots,a_n)$, $h=(h_1,\cdots,h_{n+1})$ and $\gamma=\{v_k\}_{k=1}^{n+1}\subset \mathbb{R}^n$, the notation of Lemmas \ref{MP_lem1} and \ref{MP_lem2}. If we impose that $v_k=e_k$ for $k=1,2,\cdots,n$ and  $v_{n+1}=e_s$, where $e_i=(0,0,\cdots,1^{(\text{i-th position})},0,\cdots,0)\in\mathbb{R}^n$, $i=1,\cdots,n$, then
\[
\varphi_{\gamma}(t_1,\cdots,t_{n+1})=(a+(t_1,\cdots,t_{s-1},t_s+t_{n+1},t_{s+1},\cdots,t_n),P_{a,h,\gamma}(t_1,\cdots,t_{n+1})),
\]
and
\[
P_{a,h,\gamma}(0,\cdots,0,t_s,0,\cdots,0,t_{n+1})=p_{a_s,h_s,h_{n+1}}(t_s,t_{n+1})=\sum_{k=0}^NA_k(t_s+t_{n+1})t_s^k.
\]
\end{lemma}

\noindent \textbf{Proof. } It is trivial. The result follows just by imposing the substitutions $v_k=e_k$ for $k=1,2,\cdots,n$ and $v_{n+1}=e_s$ and using the definition of 
 $p_{a_s,h_s,h_{n+1}}$ as an interpolation polynomial. {\hfill $\Box$}

\noindent \textbf{Proof of Theorem \ref{MP_grafo_vv} } The first equality from Lemma \ref{MP_lem3} implies that 
\begin{equation}\label{MP_det}
\varphi_{\gamma}'=\left [
\begin{array}{cccccccc}
1 & 0 & 0 & \cdots & 0 &\cdots & 0 & 0 \\
0 & 1 & 0 & \cdots & 0& \cdots & 0 &  0 \\
\vdots & \vdots & \vdots & \ddots & \vdots &\cdots & \vdots & \vdots \\
0 & 0 & 0 & \cdots & 1 & \cdots & 0 &  1\\
\vdots & \vdots & \vdots & \ddots & \vdots &\cdots & \vdots & \vdots \\
0 & 0 & 0 & \cdots & 0 &\cdots & 1 &  0\\
\frac{\partial P_{a,h,\gamma}}{\partial t_1} & \frac{\partial P_{a,h,\gamma}}{\partial t_2}  & \frac{\partial P_{a,h,\gamma}}{\partial t_3}  & \cdots & \frac{\partial P_{a,h,\gamma}}{\partial t_s} & \cdots & \frac{\partial P_{a,h,\gamma}}{\partial t_n}  &  \frac{\partial P_{a,h,\gamma}}{\partial t_{n+1}} 
\end{array} \right],
\end{equation}
so that, developing the determinant $\det\varphi_{\gamma}'$  by its last file, and using the notation $P=P_{a,h,\gamma}$, we get  
\begin{eqnarray*}
&\ & \xi(t_1,\cdots,t_{n+1}) := \det\varphi_{\gamma}' (t_1,\cdots,t_{n+1})\\
&\ & \ \ =  (-1)^{n+1+s}\frac{\partial P}{\partial t_{s}}(t_1,\cdots,t_{n+1})\cdot (-1)^{n-s}+\frac{\partial P}{\partial t_{n+1}}(t_1,\cdots,t_{n+1})\\
 &\ & \ \   = (\frac{\partial P}{\partial t_{n+1}} - \frac{\partial P}{\partial t_{s}})(t_1,\cdots,t_{n+1}) .
\end{eqnarray*}
Evaluating the polynomial  $\xi$ in $(0,0,\cdots,0,t_s,0,\cdots,t_{n+1})$ and using the second equality from Lemma \ref{MP_lem3}, we get
\[
\det\varphi_{\gamma}'(0,0,\cdots,0,t_s,0,\cdots,t_{n+1})= -\sum_{k=1}^NkA_k(t_s+t_{n+1})t_s^{k-1}\neq 0.
\]
Hence $\det\varphi_{\gamma}'(t_1,t_2,\cdots,t_{n+1})$ is a nonzero algebraic polynomial in the variables $t_1,\cdots,t_{n+1}$. Thus, the associated algebraic variety 
\[
Z(\det\varphi_{\gamma}')=\{(\alpha_1,\cdots,\alpha_{n+1})\in\mathbb{R}^{n+1}: \det\varphi_{\gamma}'(\alpha_1,\cdots,\alpha_{n+1})=0\}
\]
is a proper closed subset of $\mathbb{R}^{n+1}$ with empty interior. Thus,   
$\Omega=\mathbb{R}^{n+1}\setminus Z(\det\varphi_{\gamma}')$ is an unbounded open set and  the Open Mapping Theorem for differentiable functions defined on Euclidean vector spaces implies that 
$\varphi_{\gamma}(\Omega)$ is an open subset of $\mathbb{R}^{n+1}$ which is contained into $\overline{G(f)}^{\mathbb{R}^{n+1}}$, which is what we were looking for.  The part of the theorem which claims that $\varphi_{\gamma}(\Omega)$  is locally unbounded follows directly from the second equality from Lemma \ref{MP_lem3}. 

{\hfill $\Box$}

 \bibliographystyle{amsplain}


\bigskip

\footnotesize{J. M. Almira and Kh. F. Abu-Helaiel

Departamento de Matem\'{a}ticas. Universidad de Ja\'{e}n.

E.P.S. Linares,  C/Alfonso X el Sabio, 28

23700 Linares (Ja\'{e}n) Spain

Email: jmalmira@ujaen.es; kabu@ujaen.es }

\end{document}